\documentclass[a4paper,11pt,centertags,psamsfonts]{amsart}
\usepackage{amssymb}
\usepackage{times}
\usepackage[mathcal]{euscript}

\DeclareSymbolFont{operators}   {OT1}{ptmcm}{m}{n}
\DeclareSymbolFont{letters}     {OML}{ptmcm}{m}{it}
\DeclareSymbolFont{symbols}     {OMS}{pzccm}{m}{n}
\DeclareSymbolFont{largesymbols}{OMX}{psycm}{m}{n}

\DeclareMathAlphabet{\mathsf}{OT1}{phv}{m}{n}
\DeclareMathAlphabet{\mathrm}{OT1}{ptm}{m}{n}

\DeclareSymbolFont{ER}{U}{eur}{m}{n}
\DeclareSymbolFont{SY}{U}{psy}{m}{n}

\DeclareMathSymbol{,}{\mathpunct}{SY}{'054}
\DeclareMathSymbol{.}{\mathpunct}{SY}{'056}
\DeclareMathSymbol{:}{\mathpunct}{SY}{'072}

\DeclareMathSymbol{(}{\mathopen}{SY}{'050}
\DeclareMathSymbol{)}{\mathclose}{SY}{'051}

\DeclareMathSymbol{+}{\mathbin}{SY}{'053}
\DeclareMathSymbol{-}{\mathbin}{SY}{'055}
\DeclareMathSymbol{=}{\mathbin}{SY}{'075}
\DeclareMathSymbol{<}{\mathbin}{SY}{'074}
\DeclareMathSymbol{>}{\mathbin}{SY}{'076}
\DeclareMathSymbol{\leq}{\mathbin}{SY}{'243}
\DeclareMathSymbol{\geq}{\mathbin}{SY}{'263}
\DeclareMathSymbol{\nneq}{\mathbin}{SY}{'271}
\DeclareMathSymbol{\in}{\mathbin}{SY}{'316}
\DeclareMathSymbol{\nnotin}{\mathbin}{SY}{'317}
\DeclareMathSymbol{\times}{\mathbin}{SY}{'264}
\DeclareMathSymbol{\pm}{\mathbin}{SY}{'261}
\DeclareMathSymbol{\subset}{\mathbin}{SY}{'314}
\DeclareMathSymbol{\supset}{\mathbin}{SY}{'311}
\DeclareMathSymbol{\subseteq}{\mathbin}{SY}{'315}
\DeclareMathSymbol{\supseteq}{\mathbin}{SY}{'312}

\DeclareMathSymbol{/}{\mathord}{SY}{'057}
\DeclareMathSymbol{|}{\mathord}{SY}{'174}
\DeclareMathSymbol{\ast}{\mathord}{SY}{'052}
\DeclareMathSymbol{\perp}{\mathord}{SY}{'136}

\newcommand{\R}{\mathbb{R}}

\newcommand{\N}{\mathbb{N}}

\newcommand{\cI}{\mathcal{I}}

\newcommand{\cC}{\mathcal{C}}
\newcommand{\cQ}{\mathcal{Q}}
\newcommand{\cJ}{\mathcal{J}}

\newcommand{\fH}{\mathfrak{H}}

\newcommand{\sign}{\mathrm{sign}}

\DeclareMathOperator{\tr}{tr}
\DeclareMathOperator{\spec}{spec}

\newtheorem{theorem}{Theorem}[section]{\bf}{\it}
\newtheorem{proposition}[theorem]{Proposition}{\bf}{\it}
{\bf}{\it}
{\bf}{\it}
\newtheorem{lemma}{Lemma}[section]{\bf}{\it}
{\bf}{\it}
{\it}{\rm}

\newtheorem{remark}{Remark}[section]{\it}{\rm}

\newtheorem{introtheorem}{Theorem}{\bf}{\it}
\newtheorem{introcorollary}{Corollary}{\bf}{\it}

\newcommand{\hnl}{\texttt}

\title[Concavity of the Spectral Shift Function]
{Concavity of Eigenvalue Sums and the Spectral Shift Function}

\author[Vadim Kostrykin]{Vadim Kostrykin}

\address{Fraunhofer-Institut f\"{u}r Lasertechnik, Steinbachstra{\ss}e 15, D-52074 Aachen,
Germany}

\email{kostrykin@t-online.de, kostrykin@ilt.fhg.de}

\dedicatory{Dedicated to Robert Schrader on the occasion of his 60th
birthday}



\thanks{Published in J. Funct. Anal. \textbf{176} (2000), 100 -- 114.}

\begin{document}

\begin{abstract}
It is well known that the sum of negative (positive) eigenvalues of some finite
Hermitian matrix $V$ is concave (convex) with respect to $V$. Using the theory
of the spectral shift function we generalize this property to self-adjoint
operators on a separable Hilbert space with an arbitrary spectrum. More
precisely, we prove that the spectral shift function integrated with respect to
the spectral parameter from $-\infty$ to $\lambda$ (from $\lambda$ to
$+\infty$) is concave (convex) with respect to trace class perturbations. The
case of relative trace class perturbations is also considered.
\end{abstract}

\maketitle
\thispagestyle{empty}

\section{Introduction and Main Results}

Consider an arbitrary Hermitian matrix $V$. Let $\lambda_j(V)$, $j\in\N$ denote
its eigenvalues enumerated in the increasing order and repeated according to
their multiplicity. Consider the following eigenvalue sums
\begin{displaymath}
S^{(-)}_\lambda(V)=\sum_{j:\lambda_j(V)\leq\lambda}(\lambda_j(V)-\lambda),\qquad
S^{(+)}_\lambda(V)=\sum_{j:\lambda_j(V)\geq \lambda}(\lambda_j(V)-\lambda),
\end{displaymath}
which equivalently can be written in the form
\begin{eqnarray*}
S_{\lambda}^{(-)}(V)=
-\int_{-\infty}^\lambda N^{(-)}(\lambda';V)d\lambda',\qquad
S_{\lambda}^{(+)}(V)=
\int_\lambda^\infty N^{(+)}(\lambda';V)d\lambda',
\end{eqnarray*}
with $N^{(\pm)}$ being the counting functions, i.e., $N^{(\pm)}(\lambda;V)=\#\{j:\
\pm\lambda_j(V)\geq \pm\lambda\}$. By means of the min-max principle it can be easily
proved (see, e.g., \cite{Th,Lieb:Siedentop}) that $S^{(-)}(V)$ is concave and
$S^{(+)}(V)$ is convex with respect to $V$, i.e., for any Hermitian matrices $V_1$
and $V_2$ and any $\alpha\in[0,1]$
\begin{equation}\label{properties}
\pm\ S_\lambda^{(\pm)}(\alpha V_1+(1-\alpha)V_2)\ \leq\ \pm\left(
\alpha\ S_\lambda^{(\pm)}(V_1)+ (1-\alpha)\ S_\lambda^{(\pm)}(V_2)\right).
\end{equation}
These inequalities play an important role in several problems of quantum and
statistical physics (see, e.g., references cited in \cite{Lieb:Siedentop}).

In the present note we show that for a wide class of self-adjoint
operators on a separable Hilbert space $\fH$, which need not have
purely discrete spectrum, the properties \eqref{properties} are
valid for properly regularized $S_{\lambda}^{(\pm)}$. More
precisely, instead of $V$ compared to the zero operator we
consider pairs $(A_0+V, A_0)$.  For an arbitrary self-adjoint
operator $A_0$ and any self-adjoint trace class operator $V$ we
define
\begin{equation}\label{concav:1}
\zeta^{(-)}(\lambda;A_0+V,A_0)\ :=\ \int_{-\infty}^\lambda \xi(\lambda';A_0+V,A_0)d\lambda'
\end{equation}
and
\begin{equation}\label{convex:1}
\zeta^{(+)}(\lambda;A_0+V,A_0)\ :=\ \int_\lambda^\infty \xi(\lambda';A_0+V,A_0)d\lambda',
\end{equation}
where $\xi(\lambda;A_0+V,A_0)$ is the spectral shift function for the pair of
operators ($A_0+V,A_0$). Recall that for an arbitrary self-adjoint operator $A_0$ and
any self-adjoint trace class operator $V$ the spectral shift function
$\xi(\lambda;A_0+V,A_0)$ exists such that $\xi(\cdot;A_0+V,A_0)\in L^1(\R)$. Let
$F\in C^1_{\mathrm{loc}}(\R)$ be such that its derivative $F'$ belongs to the Wiener
class $W(\R)$, i.e., $F'(\lambda)$ is representable in the form
\begin{displaymath}
F'(\lambda)=\int_\R e^{-i\lambda t}d\sigma(t),
\end{displaymath}
where $\sigma(\cdot)$ is a finite complex-valued Borel measure on $\R$,
$|\sigma(\R)|<\infty$. Then
\begin{equation}\label{Krein1}
\tr\left(F(A_0+V)-F(A_0)\right)=\int_\R F'(\lambda)\xi(\lambda;A_0+V,A_0)d\lambda.
\end{equation}
This last equation may be used as a definition of the spectral shift function.
A wider class of functions for which the trace formula (\ref{Krein1}) remains
valid is discussed in \cite{Birman:Solomyak:75}. A review on the spectral shift
function is the paper by Birman and Yafaev \cite{Birman:Yafaev} (see also the
book \cite{Yafaev} and
\cite{Pushnitski:97,Gesztesy:Makarov:Naboko,Gesztesy:Makarov:99,Gesztesy:Makarov:99b}
for recent results).

In the sequel we use the notation $\cJ_p$, $p\geq 1$ for the von Neumann
- Schatten ideals of compact operators such that in particular $\cJ_1$ denotes
the set of the trace class operators (see, e.g., \cite{Gohberg:Krein}).
$\spec(A)$ denotes the spectrum of the operator $A$. $\cQ(A)$ is the domain of
the quadratic form associated with the self-adjoint operator $A$.

If in some open interval $(a,b)$ the spectrum of $A_0$ is purely discrete then
$\xi(b-0;A_0+V,A_0)-\xi(a+0;A_0+V,A_0)$ equals the difference of the total
multiplicities of the spectra of $A_0$ and $A_0+V$ lying in $(a,b)$. Thus if we
take $A_0=\lambda_+ I$ with some $\lambda_+>\sup\spec(V)$, then
$\zeta^{(-)}(\lambda;A_0+V,A_0)=S_\lambda^{(-)}(V)$ for all
$\lambda<\lambda_+$. Similarly $A_0=\lambda_- I$ with some
$\lambda_-<\inf\spec(V)$ leads to
$\zeta^{(+)}(\lambda;A_0+V,A_0)=S_\lambda^{(+)}(V)$ for all
$\lambda>\lambda_-$.

\begin{introtheorem}\label{Th1}
Let $A_0$ and $V$ be self-adjoint operators on a separable Hilbert space $\fH$,
$V\in\cJ_1$. For an arbitrary real-valued nonincreasing $f$ of bounded total
variation the functional
\begin{equation}\label{gV}
g(V)=\int_\R f(\lambda)\xi(\lambda;A_0+V,A_0)d\lambda
\end{equation}
is concave with respect to the perturbation $V$, i.e., for arbitrary $V_1,\
V_2\in\cJ_1$ the inequality
\begin{equation}\label{prop1}
g(\alpha V_1+(1-\alpha)V_2)\ \geq\ \alpha g(V_1)+(1-\alpha)g(V_2)
\end{equation}
holds for all $\alpha\in[0,1]$.
\end{introtheorem}

In particular we can take $f(\lambda)=\chi_{(-\infty,\lambda_0]}(\lambda)$, the
characteristic function of $(-\infty,\lambda_0]$ with arbitrary
$\lambda_0\in\R$, such that $g(V)=\zeta^{(-)}(\lambda_0;A_0+V,A_0)$, the
integrated spectral shift function (\ref{concav:1}). From Theorem \ref{Th1} it
follows that $\zeta^{(-)}(\lambda;A_0+V,A_0)$ is concave with respect to $V$.

It is known that
\begin{displaymath}
\int_\R\xi(\lambda; A_0+V,A_0)d\lambda=\tr V,
\end{displaymath}
which is obviously linear in $V$. Since an arbitrary nondecreasing function
$\widetilde{f}$ of bounded total variation can be represented as a difference
of a constant and a nonincreasing $f$ of bounded total variation we obtain

\begin{introcorollary}\label{Cor1}
Let $A_0$ and $V$ be as in Theorem \ref{Th1}. For an arbitrary real-valued
nondecreasing $\widetilde{f}$ of bounded total variation the functional
\begin{displaymath}
\widetilde{g}(V)=\int_\R\widetilde{f}(\lambda)\xi(\lambda;A_0+V,A_0)d\lambda
\end{displaymath}
is convex with respect to the perturbation $V$, i.e., for arbitrary $V_1,\
V_2\in\cJ_1$ the inequality
\begin{displaymath}
\widetilde{g}(\alpha V_1+(1-\alpha)V_2)\leq \alpha \widetilde{g}(V_1)+
(1-\alpha)\widetilde{g}(V_2)
\end{displaymath}
holds for all $\alpha\in[0,1]$.
\end{introcorollary}

In particular $\widetilde{f}(\lambda)=\chi_{[\lambda_0,+\infty)}(\lambda)$
satisfies the conditions of the corollary and thus
$\zeta^{(+)}(\lambda;A_0+V,A_0)$ defined by (\ref{convex:1}) is convex with
respect to $V$.

\begin{introcorollary}\label{Cor2.3}
Let the function $f$ satisfy the conditions of Theorem \ref{Th1}. Let
$V(\alpha)$ be a $\cJ_1$-valued operator family concave (in the operator sense)
with respect to $\alpha$. Then the real-valued function $\alpha\mapsto
g(V(\alpha))$ is concave. Similarly, if $f$ satisfies the conditions of
Corollary \ref{Cor1} and $V(\alpha)$ is convex, then $\alpha\mapsto
\widetilde{g}(V(\alpha))$ is also convex.
\end{introcorollary}

Theorem \ref{Th1} and Corollary \ref{Cor2.3} will be proved in Section
\ref{sec:num:2} below.

We note that a special case of this result was proved recently by Gesztesy,
Makarov, and Motovilov \cite[Corollary 1.9]{Gesztesy:Makarov:Motovilov} by
different methods.

In the article \cite{Geisler:Kostrykin:Schrader} written by the present author
in collaboration with R. Geisler and R. Schrader we have proven that the
integrated spectral shift function for the pair of Schr\"{o}dinger operators is
concave with respect to the perturbation potential. Here we will prove that
this property holds for an arbitrary pair ($A$, $A_0$) of self-adjoint
semibounded operators on a separable Hilbert space $\fH$ provided that their
difference is a relative trace class perturbation of $A_0$.

More precisely, we suppose that $A_0$ is a self-adjoint operator, semibounded from
below, and $V$ is also self-adjoint and $A_0$-compact in the form sense, i.e., for
all $a>-\inf\spec(A_0)$ the operator $(A_0+a)^{-1/2}V(A_0+a)^{-1/2}$ is compact. Then
the operator $A_V=A_0+V$, defined in the form sense, is self-adjoint with
$\cQ(A_V)=\cQ(A_0)$ and also semibounded from below. Suppose that for some $p\geq 1$
and for all sufficiently large $a$
\begin{equation}\label{incl.1}
(A_V+a)^{-p}-(A_0+a)^{-p}\in\cJ_1.
\end{equation}

If $\cI$ is an interval of the real axis such that
$\cI\supset\spec(A_V)\cup\spec(A_0)$ and for some real-valued strictly monotone
$\varphi\in C^2_{\mathrm{loc}}(\cI)$ the difference $\varphi(A_V)-\varphi(A_0)$
is trace class then the spectral shift function $\xi(\lambda; A_V, A_0)$ for
the pair of operators $(A_V, A_0)$ can be defined by means of the relation
\begin{equation}\label{inv:prinz}
\xi(\lambda; A_V, A_0) := \epsilon\ \xi(\varphi(\lambda); \varphi(A_V),\varphi(A_0)),\quad
\epsilon= \sign\ \varphi'(\lambda),
\end{equation}
which turns out to be independent of $\varphi$. Obviously, $\xi(\lambda; A_V,
A_0)$ satisfies the trace formula
\eqref{Krein1} for some class of admissible functions $F$. This construction is
known in the literature as the ``invariance principle" for the spectral shift
function (see, e.g., \cite{Birman:Yafaev,Yafaev}). Setting in \eqref{inv:prinz}
$\varphi(\lambda)=(\lambda+a)^{-p}$ we obtain
\begin{equation}\label{invariance}
\xi(\lambda;A_V,A_0)=-\xi((\lambda+a)^{-p}; (A_V+a)^{-p}, (A_0+a)^{-p}).
\end{equation}
It vanishes for all $\lambda<\inf\{\spec(A_V), \spec(A_0)\}$.

In the special case $A_0=-\Delta$ in $L^2(\R^\nu)$ and $V$ being the
multiplication operator by a real-valued measurable function $V(x)$ the
conditions above are satisfied with any $p>(\nu-1)/2$ for $\nu\geq 4$ and with
$p=1$ for $\nu\leq 3$ provided that
\begin{equation*}
\begin{array}{lcl}
V\in L^{\nu/2}(\R^\nu)\cap l^1(L^2(\R^\nu))& \mathrm{for} & \nu\geq 5,\\ V\in
L^{r}(\R^\nu)\cap l^1(L^2(\R^\nu))& \mathrm{for} & \nu=4\quad \mathrm{and\
some}\ r>2,\\ V\in L^{2}(\R^\nu)\cap L^1(\R^\nu)& \mathrm{for} & \nu=2, 3,\\
V\in L^1(\R)& \mathrm{for} & \nu=1.
\end{array}
\end{equation*}
For the definition of the Birman - Solomyak classes $l^p(L^q)$ see, e.g.,
\cite{Simon:79a}.

Let $\cC(A_0,a_0,p)$, $a_0\in\R$, $p\geq 1$ denote a set of self-adjoint
operators on the separable Hilbert space $\fH$ satisfying the following
properties:

(i) every $V\in\cC(A_0,a_0,p)$ is $A_0$-compact in the form sense;

(ii) $a_0>-\inf \spec(A_V)$ for all $V\in\cC(A_0,a_0,p)$ and  the condition
\eqref{incl.1} is satisfied for all $V\in\cC(A_0,a_0,p)$ and all $a\geq a_0$;

(iii) the set $\cC(A_0,a_0,p)$ is convex, i.e., $V_1,V_2\in\cC(A_0,a_0,p)$ implies
that $\alpha V_1+(1-\alpha)V_2\in\cC(A_0,a_0,p)$ for all $\alpha\in[0,1]$.

We will say that a set possessing these properties for some $a_0\in\R$ and
$p\geq 1$ is $A_0$-\emph{convex}. Obviously $\cC(A_0,a_0,p)$ is also $A$-convex
for any operator $A$ such that $A-A_0\in\cC(A_0,a_0,p)$.

As an example consider two self-adjoint operators $V_j$ which are $A_0$-compact
in the form sense and satisfy
\begin{equation*}
(A_0+a)^{-1/2}V_j (A_0+a)^{-p-1/2}\in\cJ_1,\qquad j=1,2
\end{equation*}
for some $a>-\inf\spec(A_0)$ and $p\geq 1$. Any operator lying in the convex
hull $\{\alpha V_1+(1-\alpha)V_2,\ \alpha\in[0,1]\}$ of $\{V_1, V_2\}$ is
obviously also $A_0$-compact in the form sense. Take $a_0>a$ such that
\begin{equation*}
\|(A_0+a_0)^{-1/2} V_j (A_0+a_0)^{-1/2}\|<1
\end{equation*}
for both $j=1,2$. By Theorem XI.12 of \cite{RS3} we obtain that the condition
\eqref{incl.1} is satisfied for all $V=\alpha V_1+(1-\alpha)V_2$ with $\alpha\in[0,1]$
and arbitrary $a\geq a_0$. Thus, the convex hull of $\{V_1, V_2\}$ is
$A_0$-convex.

\begin{introtheorem}\label{Th2}
Let $A_0$ be a self-adjoint operator semibounded from below and
$\cC(A_0,a_0,p)$ be some $A_0$-convex set. Let $q$ equal $p$ if $p=1$ and the
smallest odd integer larger than $p$ if $p>1$. Let $A_V$ with
$V\in\cC(A_0,a_0,p)$ denote the operator $A_0+V$ defined in the form sense. For
an arbitrary real-valued nonnegative nonincreasing $f$ of bounded total
variation on $[-a_0,+\infty)$ such that
\begin{displaymath}
\sup_{\lambda\in[-a_0,+\infty)}(1+|\lambda|)^{q+1}|f(\lambda)|<\infty
\end{displaymath}
the functional
\begin{displaymath}
g(V)=\int_\R f(\lambda)\xi(\lambda;A_V,A_0)d\lambda
\end{displaymath}
is concave on $\cC(A_0,a_0,p)$, i.e., for arbitrary $V_1,\ V_2\in\cC(A_0,a_0,p)$ the
inequality
\begin{equation}\label{prop2}
g\left(\alpha V_1+(1-\alpha)V_2\right)\ \geq\ \alpha g(V_1)+(1-\alpha)g(V_2)
\end{equation}
holds for all $\alpha\in[0,1]$.
\end{introtheorem}

The proof of this theorem will be given in Section \ref{sec:num:3} below.

As discussed in \cite{Geisler:Kostrykin:Schrader} (see also Proposition
\ref{vyvod2} below) the concavity (convexity) of $g(V)$ ($\widetilde{g}(V)$,
respectively) implies that $g(\alpha V)$ is subadditive and
$\widetilde{g}(\alpha V)$ is superadditive with respect to $\alpha\in\R_+$.
Subadditivity and superadditivity properties with respect to the
\emph{perturbation} (rather than with respect to the coupling constant) do not
hold generally. In the special case of the Schr\"{o}dinger operators this was
observed in
\cite{Geisler:Kostrykin:Schrader,Kostrykin:Schrader:99a,Kostrykin:Schrader:99d}.
Subadditivity and superadditivity properties of the spectral shift function
play an important role in some problems related to random Schr\"{o}dinger operators
\cite{Kostrykin:Schrader:99a,Kostrykin:Schrader:99d}. Also they allow one to
study the strong coupling limit. In particular, in Section \ref{sec:num:3} we
will prove

\begin{introcorollary}\label{cor5}
Let $A_0$ be an arbitrary self-adjoint operator and $V\geq 0$. Assume that
either

(i) $V$ is trace class

\noindent or

(ii) $A_0$ is semibounded from below, $V$ is $A_0$-compact in the form sense
and $(A_0+a)^{-1/2}V$ $\cdot(A_0+a)^{-p-1/2}\in\cJ_1$ for some $p\geq 1$ and
some $a>-\inf\spec(A_0)$.

Then for any nonincreasing function $f$ of bounded total variation, which in
the case (ii) satisfies additionally the conditions of Theorem \ref{Th2}, the
limit
\begin{equation*}
\lim_{\alpha\rightarrow\infty}\frac{1}{\alpha}\int_\R f(\lambda)
\xi(\lambda; A_0+\alpha V, A_0) d\lambda
\end{equation*}
exists and is finite.
\end{introcorollary}

For other results related to the strong coupling limit we refer to
\cite{Pushnitski:98a,Pushnitski:98b,Safronov}.

Most of the results of the present note have appeared previously in
\cite{Kostrykin} in a slightly less general form.

\section{Trace Class Perturbations}\label{sec:num:2}
\renewcommand{\theequation}{\arabic{section}.\arabic{equation}}
\setcounter{equation}{0}

The proof of Theorem \ref{Th1} relies on the following result of Birman and
Solomyak \cite{Birman:Solomyak:75}:

\begin{lemma}\label{BS}
Let $f\geq 0$ be a nonincreasing function with bounded total variation. Then
for any self-adjoint operators $A_0$ and $V$ on $\fH$, $V\in\cJ_1$

(i) the real-valued function $\alpha\mapsto\tr\left[f(A_0+\alpha V)V\right]$ is
nonincreasing, i.e., for $\alpha_1\leq\alpha_2$ the inequality
\begin{displaymath}
\tr\left[f(A_0+\alpha_1 V)V\right]\ \geq\ \tr\left[f(A_0+\alpha_2 V)V\right]
\end{displaymath}
holds,

(ii)
\begin{displaymath}
\int_\R f(\lambda)\xi(\lambda;A_0+\alpha V, A_0)d\lambda\ =\
\int_0^\alpha\tr\left[f(A_0+sV)V\right]ds.
\end{displaymath}
\end{lemma}

\vspace{2mm}

\begin{remark}
The proof in \cite{Birman:Solomyak:75} of the part (i) relies on the theory of
the double Stieltjes operator integral. An alternative proof not using this
formalism is given by Gesztesy, Makarov, and Motovilov in
\cite{Gesztesy:Makarov:Motovilov}. Part (ii) of the lemma is proven in
\cite{Birman:Solomyak:75} for the case
$f(\lambda)=\chi_{(-\infty,\lambda_0]}(\lambda)$. The present extension is
immediate. Alternative proofs of (ii) have appeared in
\cite{Simon:preprint:96,Gesztesy:Makarov:Naboko}. An operator-valued version of
this formula for sign-definite perturbations is given in
\cite{Gesztesy:Makarov:Naboko}.
\end{remark}

From Lemma \ref{BS} (i) it follows that
\begin{displaymath}
G:\ \alpha\mapsto\int_0^\alpha\tr\left[f(A_0+sV)V\right]ds
\end{displaymath}
is concave. Indeed a necessary and sufficient condition for $G(\cdot)$ to be
concave is
\begin{equation}\label{G}
2G(\alpha)-G(\alpha+h)-G(\alpha-h)\geq 0
\end{equation}
for all $\alpha\in\R$, $h\geq 0$. Since $\alpha\mapsto\tr\left[f(A_0+\alpha V)V
\right]$ is nonincreasing we have
\begin{displaymath}
\int_{\alpha}^{\alpha+h}\tr\left[f(A_0+s V)V
\right]ds-\int_{\alpha-h}^{\alpha}\tr\left[f(A_0+s V)V
\right]ds\leq 0,
\end{displaymath}
which is equivalent to (\ref{G}). Now by the claim (ii) of Lemma \ref{BS} it
follows that the functional $g(V)$ (\ref{gV}) is concave with respect to the
coupling constant.

By the chain rule for the spectral shift function (see, e.g.,
\cite{Birman:Yafaev})
\begin{displaymath}
\xi(\lambda;A_1+\alpha V,A_1)=\xi(\lambda;A_1+\alpha V, A_0)+\xi(\lambda;A_0, A_1),
\end{displaymath}
we have that
\begin{displaymath}
\int_\R f(\lambda)\xi(\lambda;A_1+\alpha V, A_0)d\lambda
\end{displaymath}
is also concave with respect to $\alpha$ for arbitrary $A_0$ and $A_1$ such
that $A_1-A_0\in\cJ_1$. Thus for arbitrary $t_1, t_2\in\R$ and arbitrary
$V\in\cJ_1$ we have
\begin{eqnarray*}
\lefteqn{\int_\R f(\lambda)\xi(\lambda;A_1+\alpha t_1 V+(1-\alpha)t_2 V, A_0)d\lambda}\\
&&\geq\ \alpha\int_\R f(\lambda)\xi(\lambda;A_1+t_1 V, A_0)d\lambda+
(1-\alpha)\int_\R f(\lambda)\xi(\lambda;A_1+t_2 V, A_0)d\lambda
\end{eqnarray*}
for all $\alpha\in[0,1]$. Taking $t_1=0$, $t_2=1$, $A_1=A_0+V_1$, and
$V=V_2-V_1$ we obtain
\begin{displaymath}
g(\alpha V_1+(1-\alpha) V_2)\ \geq\ \alpha g(V_1)+(1-\alpha) g(V_2),
\end{displaymath}
thus proving the claim of Theorem \ref{Th1}, however, under the additional
requirement that $f\geq 0$. To eliminate this requirement let us consider the
function $f_1$, which differs from $f\geq 0$ by a negative constant $c$. Since
$\tr V$ is linear in $V$, the induced functional
\begin{displaymath}
g_1(V)=\int_\R f_1(\lambda)\xi(\lambda;A_0+V, A_0)d\lambda=g(V)+c\tr V.
\end{displaymath}
is also concave in $V$. This completes the proof of Theorem \ref{Th1}.

\begin{proof}[Proof of Corollary \ref{Cor2.3}]
By Theorem \ref{Th1}
\renewcommand{\theequation}{\arabic{section}.\arabic{equation}}
\begin{equation}\label{ref1}
g(\alpha V_1+(1-\alpha)V_2)\ \geq\ \alpha g(V_1)+(1-\alpha)g(V_2)
\end{equation}
for all $\alpha\in[0,1]$. By the monotonicity of the spectral shift function with
respect to the perturbation $g(V)$ is nondecreasing with respect to $V$, i.e.,
$g(V_1)\geq g(V_2)$ for $V_1\geq V_2$.

Let now $V_1=V(\beta_1)$ and $V_2=V(\beta_2)$. By the concavity of $V(\alpha)$, i.e.,
by
\begin{displaymath}
V(\alpha\beta_1+(1-\alpha)\beta_2)\ \geq\ \alpha V(\beta_1)+(1-\alpha)V(\beta_2),
\end{displaymath}
and by the monotonicity of $g(V)$, from (\ref{ref1}) it follows that
\begin{displaymath}
g(V(\alpha\beta_1+(1-\alpha)\beta_2))\ \geq\ \alpha
g(V(\beta_1))+(1-\alpha)g(V(\beta_2)).
\end{displaymath}
The second part of the claim can be proved similarly.
\end{proof}

\section{Relative Trace Class Perturbations}\label{sec:num:3}

\renewcommand{\theequation}{\arabic{section}.\arabic{equation}}
\setcounter{equation}{0}

We turn to the case of relative trace class perturbations of $A_0$ and prove
Theorem~\ref{Th2}. The conditions of this theorem imply that
\begin{equation}\label{neu}
\sup_{\lambda\in[-a_0,+\infty)}(\lambda+a_0)^{q+1}|f(\lambda)|<\infty
\end{equation}
with $q$ being equal to $p$ if $p=1$ and to the smallest odd integer larger
than $p$ if $p>1$. Choose arbitrary $V_1, V_2\in\cC(A_0,a_0,p)$. Obviously,
$\{(1-\alpha)(V_2-V_1),\
\alpha\in[0,1]\}\subseteq\cC(A_0,a_0,p)-V_1$ and the set
$\cC_1(A_0,a_0,p):=\cC(A_0,a_0,p)-V_1$ is $A_0$-convex. Note that
$0\in\cC_1(A_0,a_0,p)$. Thus, as in the case of trace class perturbations it
suffices to prove that for any $V\in\cC_1(A_0,a_0,p)$ the function $g(\alpha
V)$ is concave with respect to $\alpha\in[0,1]$.

We start with the simplest case $p=1$ in the condition
\eqref{incl.1}. For all $a\geq a_0$ the resolvents $(A_{\alpha
V}+a)^{-1}$ and $(A_{0}+a)^{-1}$ are bounded nonnegative operators. By
assumption their difference is trace class and therefore the spectral shift
function $\xi(\lambda;A_{\alpha V}, A_0)$ can be defined by means of the
invariance principle as given by (\ref{invariance}). For all $a\geq a_0$ it
satisfies the inequality
\begin{equation}\label{ineq}
\int_\R\frac{|\xi(\lambda;A_{\alpha V}, A_0)|}{(\lambda+a)^2}d\lambda\ \leq\
\|(A_{\alpha V}+a)^{-1}-(A_0+a)^{-1}\|_{\cJ_1}.
\end{equation}

\begin{lemma}\label{Lem3.1}
Let $f(\lambda)\geq 0$ be a nonincreasing function of bounded total variation.
Then for all $a\geq a_0$
\begin{displaymath}
g_a(\alpha)=\int_\R\frac{f(\lambda)}{(\lambda+a)^2}\xi(\lambda;A_{\alpha V},
A_0)d\lambda
\end{displaymath}
is concave with respect to $\alpha$.
\end{lemma}

\begin{proof}
We change the integration variable $\lambda\rightarrow t=(\lambda+a)^{-1}$ and
use the invariance principle (\ref{invariance}) to obtain
\begin{equation}\label{gaalph}
g_a(\alpha)=-\int_0^\infty f\left(\frac{1-at}{t}\right)\xi(t;(A_{\alpha
V}+a)^{-1},(A_0+a)^{-1})dt.
\end{equation}
It is easy to see that $f((1-at)/t)$ is nondecreasing with respect to $t$. It is well
known (see, e.g., \cite{Kraus,Bendat:Sherman}, \cite[Proposition 1.3.11]{Pedersen})
that the function $x\mapsto x^{-1}$ is concave on the set of invertible positive
operators, i.e., for arbitrary invertible positive operators $X$ and $Y$ the
inequality $(\beta X+(1-\beta) Y)^{-1}\leq
\beta X^{-1}+(1-\beta) Y^{-1}$ holds in the operator sense  for all $\beta\in[0,1]$.
Taking $X=A_{\alpha_1 V}+a$ and $Y=A_{\alpha_2 V}+a$ with an arbitrary $a\geq
a_0$ and using the fact that $\beta X+(1-\beta)Y=A_{(\beta\alpha_1
+(1-\beta)\alpha_2)V}+a$ we obtain that
\begin{equation*}
(A_{(\beta\alpha_1 +(1-\beta)\alpha_2)V}+a)^{-1}\ \leq\
\beta(A_{\alpha_1 V}+a)^{-1}+ (1-\beta)(A_{\alpha_2 V}+a)^{-1}
\end{equation*}
for all $\beta\in[0,1]$, i.e., the operator $(A_{\alpha V}+a)^{-1}$ is convex with
respect to $\alpha\in\R$. Therefore by Corollary \ref{Cor2.3} the integral in
(\ref{gaalph}) is convex with respect to $\alpha$ and thus $g_a(\alpha)$ is concave.
\end{proof}

From \eqref{neu} it follows that the function $f$ satisfies the condition
\begin{displaymath}
\sup_{\lambda\in[-a_0,+\infty)}(\lambda+a_0)^2|f(\lambda)|<\infty.
\end{displaymath}
Since $\xi(\lambda;A_{\alpha V}, A_0)=0$ for all $\lambda\leq -a_0$ we may
suppose that $\lambda\geq -a_0$. Thus for all $a\geq 2 a_0$ we have
$a^2(a+\lambda)^{-2}\leq 4$. Obviously,
\begin{eqnarray*}
\lefteqn{\left|\frac{a^2}{(\lambda+a)^2} f(\lambda) \xi(\lambda; A_{\alpha V}, A_0)
\right|}\\
&\leq & 4 \frac{|\xi(\lambda; A_{\alpha V}, A_0)|}{(\lambda+a_0)^2}\
\sup_{\lambda\in[-a_0,+\infty)}(\lambda+a_0)^2 |f(\lambda)|.
\end{eqnarray*}
Therefore, by (\ref{ineq}) and by the Lebesgue dominated convergence theorem we
have
\begin{displaymath}
\lim_{a\rightarrow+\infty}a^2g_a(\alpha)=\int_\R f(\lambda)\xi(\lambda;A_{\alpha V}, A_0)
d\lambda.
\end{displaymath}
From Lemma \ref{Lem3.1} it follows now that the integral on the r.h.s.\ is
concave with respect to $\alpha$. As noted above the concavity with respect to
the coupling constant implies the concavity with respect to the perturbation.
This remark completes the proof of Theorem \ref{Th2} in the case $p=1$.

We turn to the case $p>1$ in the condition \eqref{incl.1} and note that the
operator $(A_{\alpha V}+a_0)^{-p}$ is neither convex nor concave with respect
to $\alpha$ \cite{Kraus,Bendat:Sherman}. To treat this case we need the
following

\begin{lemma}$\mathrm{(}$\cite[Theorem 1]{Koplienko}$\mathrm{;}$ \cite[Theorem 8.10.4]{Yafaev}$\mathrm{)}$\label{xxx}
Assume that $A^p-A_0^p\in\cJ_1$ for some $p>1$. Then $A-A_0\in\cJ_q$ for any
$q>p$. Let $\{P_n\}_{n\in\N}$ be a strictly monotone family of finite
dimensional orthogonal projections converging strongly to the identity operator
$I$. Then
\begin{equation*}
\xi(\lambda; A, A_0) = \lim_{n\rightarrow\infty} \xi(\lambda; A_0+P_n(A-A_0)P_n, A_0)
\end{equation*}
in $L^1(\R;\lambda^{q-1}d\lambda)$, where $q$ is the smallest odd integer
greater than $p$.
\end{lemma}

The following lemma generalizes Lemma \ref{Lem3.1} to the case $p>1$:

\begin{lemma}\label{relat:neu}
Let $f(\lambda)\geq 0$ be a nonincreasing function of bounded total variation.
Then for all $a\geq a_0$
\begin{displaymath}
g_a(\alpha)=\int_\R\frac{f(\lambda)}{(\lambda+a)^{q+1}}\xi(\lambda;A_{\alpha
V}, A_0)d\lambda
\end{displaymath}
is concave with respect to $\alpha$.
\end{lemma}

\begin{proof}
By the invariance principle the spectral shift function $\xi(\lambda;A_V, A_0)$
can be represented in the form
\begin{equation*}
-\xi((\lambda+a)^{-1};(A_V+a)^{-1},(A_0+a)^{-1}).
\end{equation*}
We introduce the operator $W(\alpha)=(A_{\alpha V}+a)^{-1}-(A_0+a)^{-1}$. Let
$\{P_n\}_{n\in\N}$ be a family of finite dimensional orthogonal projections as
in Lemma \ref{xxx}. Consider $W_n(\alpha)=P_n W(\alpha) P_n\in\cJ_1$ and define
\begin{eqnarray*}
\lefteqn{g_a^{(n)}(\alpha)}\\
&=&-\int_\R\frac{f(\lambda)}{(\lambda+a)^{q+1}}\xi((\lambda+a)^{-1};
(A_0+a)^{-1}+W_n(\alpha), (A_0+a)^{-1})d\lambda\\ &=& -\int_0^\infty
f\left(\frac{1-at}{t}\right)t^{q-1}\xi(t;(A_0+a)^{-1}+W_n(\alpha),(A_0+a)^{-1})dt
\end{eqnarray*}
with $q$ being defined as in Theorem \ref{Th2}. Recall that $W(\alpha)$ is
convex with respect to $\alpha$ and therefore $P_nW(\alpha)P_n$ is also convex.
Thus by Corollary \ref{Cor2.3} the function $g_a^{(n)}(\alpha)$ is concave for
every $n\in\N$. To prove that
\renewcommand{\theequation}{\arabic{section}.\arabic{equation}}
\begin{equation}\label{Ziel}
\lim_{n\rightarrow\infty}g_a^{(n)}(\alpha)=g_a(\alpha)
\end{equation}
we estimate as follows
\begin{eqnarray*}
\lefteqn{|g_a(\alpha)-g_a^{(n)}(\alpha)|}\\ &\leq &
\sup_{\lambda\in[-a_0,+\infty)} |f(\lambda)|
\int_0^\infty t^{q-1} \Big|\xi(t;(A_0+a)^{-1}+W(\alpha),(A_0+a)^{-1})\\
&&-\xi(t;(A_0+a)^{-1}+W_n(\alpha), (A_0+a)^{-1})\Big| dt.
\end{eqnarray*}
By Lemma \ref{xxx} the r.h.s.\ tends to zero thus proving \eqref{Ziel} and
completing the proof of the lemma.
\end{proof}

To complete the proof of Theorem \ref{Th2} as in the case $p=1$ we consider the
limit $a\rightarrow+\infty$ of $a^{q+1}g_a(\alpha)$. By the inequality
\begin{eqnarray*}
\int_\R \frac{|\xi(\lambda; A_{\alpha V}, A_0)|}{(\lambda+a)^{q+1}}d\lambda
&\leq &\int_\R \frac{|\xi(\lambda; A_{\alpha V},
A_0)|}{(\lambda+a)^{p+1}}d\lambda\\ &\leq & \left\|(A_{\alpha
V}+a)^{-p}-(A_0+a)^{-p}\right\|_{\cJ_1}
\end{eqnarray*}
valid for all $a\geq a_0+1$ and again by the Lebesgue dominated convergence
theorem we obtain
\begin{equation*}
\lim_{a\rightarrow+\infty} a^{q+1} g_{a}(\alpha)=\int_\R f(\lambda) \xi(\lambda; A_{\alpha V},
A_0) d\lambda.
\end{equation*}
Now from Lemma \ref{relat:neu} it follows that the integral on the r.h.s.\ is
concave with respect to $\alpha$. This completes the proof of Theorem
\ref{Th2}.

We turn to the proof of Corollary \ref{cor5}.

\begin{proposition}\label{vyvod2}
Under the assumptions of Corollary \ref{cor5} (but without the restriction
$V\geq 0$) the functional $g(V)$ is subadditive in the coupling constant in the
sense that for arbitrary $\alpha_1,\alpha_2\geq 0$
\renewcommand{\theequation}{\arabic{section}.\arabic{equation}}
\begin{equation}\label{start}
g((\alpha_1+\alpha_2)V)\ \leq\ g(\alpha_1 V)+g(\alpha_2 V).
\end{equation}
Moreover, for arbitrary $\alpha_1,\alpha_2\geq 0$ the inequality
\begin{equation}\label{vyvod2:eq1}
g((\alpha_1-\alpha_2)V)\ \geq\ g(\alpha_1 V)+g(-\alpha_2 V)
\end{equation}
holds.
\end{proposition}

\begin{proof}
The assumption (ii) of the Corollary \ref{cor5} and the proof of Theorem XI.12
in \cite{RS3} imply that for an arbitrary finite interval $[a,b]\subset\R$
there is finite $a_0\in\R$ such that $(A_{\alpha
V}+a)^{-p}-(A_0+a)^{-p}\in\cJ_1$ for all $a\geq a_0$. Thus we may set
$\cC(A_0,a_0,p)=\{\alpha V,\ \alpha\in[a,b]\}$. By Theorem \ref{Th2} we obtain
that $g(\alpha V)$ is concave with respect to $\alpha\in[a,b]$. Since $a$ and
$b$ are arbitrary the function $g(\alpha V)$ is concave on $\R$. In the case of
assumption (i) the concavity of $g(\alpha V)$ for all $\alpha\in\R$ is
guaranteed directly by Theorem \ref{Th1}.

Recall that the necessary and sufficient condition \cite[Theorem 6.2.4]{Hille}
for a measurable concave function $\phi(\alpha)$ to be subadditive on $\R_+$ is
that $\phi(+0)\geq 0$. This proves \eqref{start}. To prove \eqref{vyvod2:eq1}
we use the fact (see, e.g., \cite[Theorem 110]{Hardy:Littlewood:Polya}) that
any continuous concave function $\phi(x)$ satisfies the inequality
\renewcommand{\theequation}{\arabic{section}.\arabic{equation}}
\begin{equation}\label{vyvod2:eq2}
\phi(x-h')+\phi(x+h')\ \geq\ \phi(x-h)+\phi(x+h)
\end{equation}
provided that $|h|\geq|h'|$. We set $x=(\alpha_1-\alpha_2)/2$,
$h'=(\alpha_2-\alpha_1)/2$, $h=(\alpha_2+\alpha_1)/2$ and apply the inequality
\eqref{vyvod2:eq2} to the function $g(\alpha V)$. Since $g(0)=0$ we arrive at the
claim \eqref{vyvod2:eq1}.
\end{proof}

\begin{proof}[Proof of Corollary \ref{cor5}]
Let $\gamma=\inf_{\alpha>0}\alpha^{-1}\phi(\alpha)$. Recall (see, e.g.,
\cite[Theorem 6.6.1]{Hille}) that if $\phi(\alpha)$ is a measurable subadditive
function, which is finite for all finite $\alpha$, then $-\infty\leq \gamma
<\infty$ and
\begin{equation*}
\lim_{\alpha\rightarrow+\infty}\frac{\phi(\alpha)}{\alpha}=\gamma.
\end{equation*}
We take $\phi(\alpha)=g(\alpha V)$. By Proposition \ref{vyvod2} it is
subadditive on $\R_+$. By the monotonicity propertiy of the spectral shift
function the condition $V\geq 0$ implies that $\phi(\alpha)\geq 0$ for all
$\alpha\in\R_+$. Therefore $\gamma\geq 0$, thus proving the corollary.
\end{proof}

\section*{acknowledgment}
The author is indebted to V. Enss and R. Schrader for useful remarks.
Stimulating communications with F. Gesztesy and K. Makarov are also gratefully
acknowledged.


\end{document}